\documentclass{amsart}
\usepackage{amssymb,amscd,enumerate,xspace,amsthm}
\usepackage[all,ps,cmtip]{xy}

\theoremstyle{plain}
\newtheorem{theoreme}{Theorem}[section]

\newtheorem{corollaire}[theoreme]{Corollary}
\newtheorem{lemme}[theoreme]{Lemma}

\theoremstyle{remark}
\newtheorem{Remarques}[theoreme]{Remarks}
\newtheorem{remarque}[theoreme]{Remark}
\newtheorem*{remarque*}{Remark}

\DeclareMathOperator{\vol}{vol}
\DeclareMathOperator{\DR}{DR}
\DeclareMathOperator{\reg}{(r)}
\DeclareMathOperator{\Reg}{reg}
\DeclareMathOperator{\ir}{(ir)}
\let\el\prime

\DeclareMathOperator{\Db}{\mathfrak{Db}}
\DeclareMathOperator{\sld}{\mathfrak{sl}_2}
\newcommand{\oxdh}{\cO_{\wh{X|D}}}
\def\mod{\mathop{\rm mod}\nolimits}
\newcommand{\Crochet}{\langle\,,\,\rangle}

\let\dpl\displaystyle
\let\wh\widehat
\let\wt\widetilde
\let\ov\overline
\let\emptyset\varnothing
\let\ldots\dots

\newcommand{\bbullet}{{\scriptscriptstyle\bullet}}

\newcommand{\ssi}{\Leftrightarrow}

\newcommand{\qqbox}[1]{\quad\hbox{#1}\quad}

\newcommand{\ccup}{\mathop\cup\limits}

\newcommand{\module}[1]{\left\vert#1\right\vert}
\newcommand{\norme}[1]{\left\Vert#1\right\Vert}

\newcommand{\defin}{\stackrel{\mathrm{def}}{=}}

\DeclareMathOperator{\id}{Id}
\DeclareMathOperator{\im}{Im}
\newcommand{\isom}{\stackrel{\sim}{\longrightarrow}}

\renewcommand{\ker}{\mathop{\rm Ker}\nolimits}

\newcommand{\reel}{\mathop{\mathrm{Re}}\nolimits}
\newcommand{\res}{\mathop{\mathrm{Res}}\nolimits}

\newcommand{\cHom}{\mathop{\cH om}\nolimits}

\DeclareMathOperator{\tr}{tr}

\newcommand{\loccit}{\emph{loc\ptbl cit}}
\newcommand{\eg}{{\it e.g}}
\newcommand{\cf}{\emph{cf}}
\newcommand{\etc}{{\it etc}}
\newcommand{\ie}{{\it i.e}}
\newcommand{\resp}{{\it resp}}

\newcommand{\T}{\S\kern .15em }
\newcommand{\ptbl}{.\kern .15em }

\newcommand{\lefpar}{\left(}
\newcommand{\rigpar}{\right)}
\newcommand{\lefcro}{\left[}
\newcommand{\rigcro}{\right]}

\newcommand{\lcr}{\left[\!\left[}
\newcommand{\rcr}{\right]\!\right]}

\newdimen\lengtharrow
 \newbox\exponantbox \newbox\indicebox
\def\dimmax#1#2{\ifdim#1<#2 #2\else #1\fi}

\def\arrowr#1#2%
    {\setbox\exponantbox=\hbox{$#2$}
      \setbox\indicebox=\hbox{$#1$}
       \lengtharrow=\dimmax{\wd\indicebox}{\wd\exponantbox}
        \ifdim \lengtharrow=0pt
                                                                                                                        \mathrel{\smash{\mathop{\hbox to 10mm{\rightarrowfill}}}}
         \else \ifdim\lengtharrow<6mm
                \lengtharrow=10mm
                 \else \advance\lengtharrow by 4mm
                  \fi
                                        \mathrel{\smash{\mathop{\hbox to\lengtharrow{\rightarrowfill}}\limits%
                         ^{\styleflechehoriz#2}_{\styleflechehoriz#1}}}
                                                                        \fi}

\def\arrowl#1#2%
    {\setbox\exponantbox=\hbox{$#2$}
      \setbox\indicebox=\hbox{$#1$}
       \lengtharrow=\dimmax{\wd\indicebox}{\wd\exponantbox}
        \ifdim \lengtharrow=0pt
                                                                                                                        \mathrel{\smash{\mathop{\hbox to 10mm{\leftarrowfill}}}}
         \else \ifdim\lengtharrow<6mm
                \lengtharrow=10mm
                 \else \advance\lengtharrow by 4mm
                  \fi
                                        \mathrel{\smash{\mathop{\hbox to\lengtharrow{\leftarrowfill}}\limits%
                         ^{\styleflechehoriz#2}_{\styleflechehoriz#1}}}
                                                                        \fi}

\def\hookarrowr#1#2{\lhook\joinrel\arrowr{#1}{#2}}

\def\dasharrowr#1#2{\mapstochar\arrowr{#1}{#2}}

\def\MRE#1{\arrowr{}{#1}}

\def\HMRE#1{\hookarrowr{}{#1}}

\def\DMRE#1{\dasharrowr{}{#1}}

\let\styleflechehoriz=\textstyle

\def\CC{\mathbf{C}}

\def\NN{\mathbf{N}}

\def\RR{\mathbf{R}}

\def\ZZ{\mathbf{Z}}

\def\cA{\mathcal{A}}

\def\cC{\mathcal{C}}
\def\cD{\mathcal{D}}
\def\cE{\mathcal{E}}

\def\cH{\mathcal{H}}

\def\cL{\mathcal{L}}
\def\cM{\mathcal{M}}
\def\cN{\mathcal{N}}
\def\cO{\mathcal{O}}

\def\cR{\mathcal{R}}

\def\cT{\mathcal{T}}

\def\varepsilong{\boldsymbol{\varepsilon}}

\def\bme{\boldsymbol{e}}

\def\bH{\boldsymbol{H}}

\def\bR{\boldsymbol{R}}

\DeclareMathAlphabet{\mathcalmaigre}{U}{eus}{m}{n}

\def\ccV{\mathcalmaigre{V}}

\DeclareMathAlphabet{\mathcalgras}{U}{eus}{b}{n}

\begin{document}

\title{Harmonic metrics and~connections~with~irregular~singularities}
\author{Claude Sabbah}
\address{UMR 7640 du CNRS\\
Centre de Math{\'e}matiques\\
{\'E}cole polytechnique\\
F--91128 Palaiseau cedex\\
France}
\email{sabbah@math.polytechnique.fr}
\urladdr{http://www.math.polytechnique.fr/cmat/sabbah/sabbah.html}
\begin{abstract}
We identify the holomorphic de~Rham complex of the minimal extension of a meromorphic vector bundle with connexion on a compact Riemann surface $X$ with the $L^2$ complex relative to a suitable metric on the bundle and a complete metric on the punctured Riemann surface. Applying results of C.~Simpson, we show the existence of a harmonic metric on this vector bundle, giving the same $L^2$ complex. As a consequence, we obtain a Hard Lefschetz-type theorem.
\end{abstract}
\subjclass{32S40, 32S60, 32L10, 35A20, 35A27}
\keywords{Harmonic metric, irregular singularity, meromorphic connection, Poincar{\'e} lemma}
\maketitle

\def\baselinestretch{1.1}\normalfont
\section{Statement of the results}\label{sec:statement}
Let $X$ be a compact Riemann surface, $D\subset X$ be a finite set of points and denote by $j$ the open inclusion $X^*\defin X-D\hookrightarrow X$. Let $\cM$ be a locally free $\cO_X[*D]$-module of finite rank $d$, equipped with a connection $\nabla:\cM\rightarrow \cM\otimes_{\cO_X}\Omega_{X}^{1}$ which may have regular or irregular singularities at each point of $D$. Therefore, $\cM$ is also a holonomic module on the ring $\cD_X$ of holomorphic differential operators on $X$. We call such a $\cD_X$-module a {\em meromorphic connection} for short. There exists a unique holonomic $\cD_X$-submodule $\cM_{\min}\subset\cM$ satisfying
\begin{enumerate}
\item
$\cO_X[*D]\otimes_{\cO_X}\cM_{\min}=\cM$,
\item
$\cM_{\min}$ has no quotient supported on a subset of $D$.
\end{enumerate}
One says that $\cM_{\min}$ is the {\em minimal extension} of $\cM$ along $D$. If $\cM^*$ denotes the dual $\cD_X$-module, we have an exact sequence
$$
0\longrightarrow \cM_{\min}\longrightarrow \cM\longrightarrow \lefcro \cH_{[D]}^{0}(\cM^*)\rigcro^*\longrightarrow 0
$$
where $\cH_{[D]}^{0}(\cM^*)$ denotes the torsion of $\cM^*$ supported on $D$.

Denote by $M$ the local system of horizontal sections of $\cM_{|X^{*}}$ and denote by $\DR(\cM)$ the de~Rham complex $(\Omega_{X}^{\bbullet}\otimes_{\cO_X}\cM,\nabla)$. If $\cM$ is regular at each point of $D$, we have $\DR(\cM_{\min})=j_*M$. In general, however, $\DR(\cM_{\min})$ cannot be computed in terms of $M$ only.

When $\nabla$ has only \emph{regular singularities}, it follows from \cite{Simpson90} that, when the meromorphic bundle with connection $(\cM,\nabla)$ is \emph{irreducible}, \ie.\ in this case when the local system $M$ is so, there exists on $\cM_{|X^*}$ a harmonic metric having a moderate behaviour near each point of $D$. Then, Zucker's arguments in \cite{Zucker79} show that the $L^2$ complex of this bundle (associated with this harmonic metric and with a metric on $X^*$ locally equivalent, near each point of $D$, to the Poincar{\'e} metric on the punctured disc) is isomorphic to the de~Rham complex $\DR(\cM_{\min})=j_*M$, and its cohomology can be computed with $L^2$ harmonic sections. It follows in particular that the class in $H^2(X,\CC)$ of a metric on $X$ (and hence, if $X$ is connected, any nonzero element in $H^2(X,\CC)=\CC$) induces a Lefschetz isomorphism
\[
\ell:\bH^0(X,\DR(\cM_{\min}))\isom \bH^2(X,\DR(\cM_{\min})).
\]
The situation is thus analogous to that of a variation of polarized Hodge structure considered by Zucker, without the Hodge decomposition however.

\medskip
In this paper, we show that the same results hold without any regularity condition on the connection. The results below are a step towards a general conjecture made by M. Kashiwara \cite{Kashiwara98} concerning the behaviour of semi-simple holonomic $\cD$-modules. The problem we solve was raised in \cite{Simpson90}.

First, we construct in \T\ref{sec:construction} a hermitian metric $k$ on the flat bundle $\cC^\infty_{X^*}\otimes_{\cO_{X^*}}\cM_{|X^*}$, satisfying the following properties:
\begin{enumerate}[({1}.1)]
\item
the $k$-norm of any local section of $\cM$ has moderate growth near each point of $D$;
\item
the norm of the curvature of $k$, computed with a $C^\infty$ metric on $X$, is ``acceptable'' in the neighbourhood of $D$ (in the sense of \cite[prop\ptbl3.1]{Simpson90});
\item
the $k$-norm of the pseudo-curvature (in the sense of \cite{Simpson90}) of $(\cM_{|X^*},\nabla,k)$ is $L^p$ with some $p>1$;
\item
the degree of the bundle $\cM_{|X^*}$, computed with the metric $k$, is zero.
\end{enumerate}

We then obtain, as a consequence of the fundamental results of Simpson \cite{Simpson88,Simpson90}:

\begin{theoreme}\label{th:harmonique}
Let $(\cM,\nabla)$ be a meromorphic bundle with connection on $X$, having poles on $D$. If $(\cM,\nabla)$ is irreducible, there exists a harmonic metric $h$ (in the sense of \cite{Simpson90}) on $\cM_{|X^*}$ such that $h$ and $k$ are mutually bounded.
\end{theoreme}

Notice that, when $\nabla$ is irregular, $(\cM,\nabla)$ may irreducible without $M$ being so.

\begin{proof}
Indeed, property (1.2) shows that the subsheaf $\wt\cM$ of $j_*\cM_{|X^*}$ of sections having a norm with moderate growth is a meromorphic bundle (\cf.\ \cite{C-G75,Simpson88}). This meromorphic bundle contains $\cM$ because of (1.1). As both bundles $\cM$ and $\wt\cM$ coincide on $X^*$, they are equal.

If $(\cM,\nabla)$ is irreducible, the object $(\cM_{|X^*},\nabla,k)$ is stable: if $\cN^*$ is a proper subbundle of $\cM_{|X^*}$ invariant by $\nabla$, Lemma 6.2 in \cite{Simpson90} shows that
\begin{itemize}
\item
either $\deg(\cN^*,k)=-\infty$ and the existence of $\cN^*$ does not contradict stability,
\item
or $\deg(\cN^*,k)$ is finite, and then $\cN^*$ can be extended as a meromorphic subbundle of $\wt\cM=\cM$; this is impossible as $(\cM,\nabla)$ is irreducible.
\end{itemize}

The object $(\cM_{|X^*},\nabla,k)$ is thus stable and has degree $0$, because of (1.4). We may apply part 2 of theorem 6 in \cite{Simpson90} to get the harmonic metric satisfying the desired properties.
\end{proof}

Fix also a complete metric on $X^*$ which, at each point of $D$, is locally equivalent to the Poincar{\'e} metric on the punctured disc $\Delta^*$. This metric and the metric $k$ on $\cM_{|X^*}$ allow us to define a $L^2$ complex on $X$, denoted by $\cL_{(2)}^\bbullet(\cM)$ (no irreducibility assumption is needed here). This complex comes equipped with a natural morphism in the complex of currents $\Db_{X}^{\bbullet}\otimes_{\cO_X}\cM$, which is isomorphic to the de~Rham complex $\DR(\cM)$.

\begin{theoreme}[$L^2$ Poincar{\'e} lemma]\label{th:poincare}
The complex $\cL_{(2)}^\bbullet(\cM)$ is isomorphic to $\DR(\cM_{\min})$.
\end{theoreme}

\begin{Remarques}\label{rem:zucker}
\begin{enumerate}\nopagebreak
\item
The de~Rham complex is centered here in the usual way, namely
\begin{eqnarray*}
\DR(\cM_{\min})&=&(\Omega_X^\bbullet\otimes \cM_{\min},\nabla).
\end{eqnarray*}
It is known to be constructible. Hence the theorem implies the constructibility of the complex $\cL_{(2)}^\bbullet(\cM)$ and the finiteness of its global cohomology.
\item
When $(\cM,\nabla)$ has only regular singularities, theorem \ref{th:poincare} is due to S. Zucker \cite{Zucker79} (th\ptbl6.2 and prop\ptbl11.3, which is also valid for $\alpha\in\CC$). In the proof, we will use this result.
\end{enumerate}
\end{Remarques}

\begin{corollaire}
Assume that $(\cM,\nabla)$ is irreducible and (for simplicity) that
$X$ is connected. Then the Lefschetz operator
\[
\ell:\bH^0(X,\DR(\cM_{\min}))\longrightarrow \bH^2(X,\DR(\cM_{\min}))
\]
induced by the cup product with any nonzero element of $H^2(X,\CC)=\CC$ is an isomorphism.
\end{corollaire}

\begin{proof}
Let $\cL_{(2)}^{\bbullet}(\cM)$ be the $L^2$ complex associated with the connection $\nabla$ and the harmonic metric $h$ given by theorem \ref{th:harmonique} (the complete metric on $X^*$ being fixed as above). This complex is equal to the $L^2$ complex associated with $\nabla$ and $k$, as $h$ and $k$ are mutually bounded, hence is isomorphic to $\DR(\cM_{\min})$, according to theorem \ref{th:poincare}. Let $\omega$ be the K{\"a}hler form associated with the fixed metric on $X^*$ and $\omega_0$ be the form associated with a $C^\infty$ metric on $X$. Standard arguments show that the previous results imply that the operator $L=\omega\wedge{}$ induces an isomorphism
\begin{eqnarray*}
\bH^0(X^*,\cL_{(2)}^{\bbullet}(\cM))&\isom&\bH^2(X^*,\cL_{(2)}^{\bbullet}(\cM)).
\end{eqnarray*}
The argument of \cite[Lemma 6.4.1]{K-K87} shows that, after the identification $\cL_{(2)}^{\bbullet}(\cM)\simeq\DR(\cM_{\min})$, this morphism corresponds to
\[
\ell:\bH^0(X,\DR(\cM_{\min}))\isom \bH^2(X,\DR(\cM_{\min}))
\]
induced by the cup product with $\omega_0$. This proves the corollary.
\end{proof}

\section{Preliminaries}\label{sec:prelim}
\subsection{}\label{subsec:rappels}
Recall first the formal structure of meromorphic bundles with connection. Let $\Delta$ be the disc $\{z\in\CC\mid\module{z}<r_0\}$. It will be convenient to assume that $r_0<1$. Let $\varphi\in\cO(\Delta)[1/z]$ and denote by $\cE^\varphi$ the free $\cO_\Delta[1/z]$-module of rank one equipped with the connection defined by $\partial _z(1)=\varphi(z)$. We say that a meromorphic connection $(\cM,\nabla)$ on $\Delta$ with pole at $0$ only is \emph{elementary} if it is isomorphic to a direct sum of meromorphic connections $\cE^\varphi\otimes_{\cO_\Delta}\cR_\varphi$ where $\cR_\varphi$ are regular at the origin.

We say that a meromorphic connection $(\cM,\nabla)$ on $\Delta$ has a \emph{formal decomposition} at the origin if there exists an elementary meromorphic connection $(\cM^{\el},\nabla^{\el})$ and an isomorphism $(\wh\cM,\wh\nabla)\simeq(\wh\cM^{\el},\wh\nabla^{\el})$, where we put, for a meromorphic connection $(\cN,\nabla)$, $\wh\cN\defin\CC\lcr z\rcr\otimes_{\cO_\Delta}\cN$. The theorem of Turrittin (see \eg. \cite{Malgrange91,Wasow65,Sibuya90bis}) asserts that any germ $(\cM,\nabla)$ has a formal decomposition, possibly after a ramification: there exists a cyclic covering
\begin{eqnarray*}
\Delta_q&\MRE{\pi}&\Delta\\
t&\longmapsto&z=t^q
\end{eqnarray*}
for some $q\in\NN^*$ such that the inverse image $\pi^*\cM$ equipped with its natural connection, that we denote by $\pi^+(\cM,\nabla)$, has a formal decomposition at the origin. There exist (see \cite[chap\ptbl III, th\ptbl (2.3)]{Malgrange91}) formal germs of regular connection $(\wh\cM^{\reg},\wh\nabla)$ and purely irregular connection $(\wh\cM^{\ir},\wh\nabla)$ such that
\begin{eqnarray*}
(\wh\cM,\wh\nabla)&=&(\wh\cM^{\reg},\wh\nabla)\oplus (\wh\cM^{\ir},\wh\nabla).
\end{eqnarray*}

One may rephrase this a little differently:

\begin{lemme}\label{lem:modele}
Let $(\cM,\nabla)$ be a meromorphic bundle with connection on $\Delta$, with pole at $0$ only. There exist
\begin{enumerate}
\item
a meromorphic bundle with connection $(\cM^{\el},\nabla^{\el})$ and an integer $q$ such that $\pi^+(\cM^{\el},\nabla^{\el})$ is elementary,
\item
an isomorphism $\wh\lambda:(\wh\cM^{\el},\wh\nabla^{\el})\isom(\wh\cM,\wh\nabla)$.
\end{enumerate}
\end{lemme}

\begin{proof}
Let $\pi:t\mapsto t^q=z$ be such that $\pi^+(\cM,\nabla)$ has a formal decomposition. Let $\zeta=e^{2i\pi/q}$ and $\sigma:t\mapsto\zeta t$. As $\pi=\pi\circ\sigma$, we have an isomorphism $\lambda_\sigma:\pi^+(\cM,\nabla)\isom\sigma^+(\pi^+(\cM,\nabla))$.

Let $(\cN^{\el},\nabla^{\el})$ be an elementary meromorphic connection in the $t$-variable equipped with a formal isomorphism $\wh\lambda_\pi:(\wh\cN^{\el},\wh\nabla^{\el})\isom\pi^+(\wh\cM,\wh\nabla)$, given by Turrittin's theorem. Put
\[
\wh{\lambda^{\el}_\sigma}=\sigma^+(\lambda_\pi)^{-1}\circ\wh{\lambda_\sigma}\circ\wh\lambda_\pi:(\wh\cN^{\el},\wh\nabla^{\el})\isom\sigma^+(\wh\cN^{\el},\wh\nabla^{\el}).
\]
As $(\cN^{\el},\nabla^{\el})$ and $\sigma^+(\cN^{\el},\nabla^{\el})$ are elementary, this isomorphism comes from an isomorphism $\lambda^{\el}_\sigma:(\cN^{\el},\nabla^{\el})\isom\sigma^+(\cN^{\el},\nabla^{\el})$. Hence, there exists a meromorphic bundle with connection $(\cM^{\el},\nabla^{\el})$ in the $z$-variable such that $(\cN^{\el},\nabla^{\el})=\pi^+(\cM^{\el},\nabla^{\el})$: it is the invariant part of the meromorphic bundle with connection $\pi_*\cN^{\el}$ with respect to the automorphism induced by $\lambda^{\el}_\sigma$.

Notice that, $\wh\lambda_\pi$ being compatible with $\wh{\lambda_\sigma}$ and $\wh\lambda^{\el}_\sigma$ by definition, it defines a formal isomorphism $\wh\lambda:(\wh\cM^{\el},\wh\nabla^{\el})\isom(\wh\cM,\wh\nabla)$.
\end{proof}

\subsection{}
Let us now be more explicit concerning $\DR(\cM_{\min})$. Denote by $\oxdh$ the formal completion of $\cO_X$ along $D$. This is a sheaf on $X$ supported on $D$, which fiber at each point of $D$ is isomorphic to $\CC\lcr z\rcr$. Consider the formalized bundle $\wh\cM\defin\oxdh\otimes_{\cO_X}\cM$. As in the local case above, it can be decomposed as the direct sum $\wh\cM=\wh\cM^{\reg}\oplus \wh\cM^{\ir}$ of its regular component and its purely irregular component.

Notice that the notion of minimal extension is also well-defined for $\oxdh[*D]$-modules with connections and that $\wh{\cM_{\min}}=\wh\cM_{\min}$: indeed, tensorizing with $\oxdh$ commutes with duality of $\cD$-modules and preserves $\cO_X$ or $\oxdh$-torsion. Denote by $\cT_\cM$ (\resp. $\cT_{\wh\cM}$) the dual of the torsion part of $\cM^*$ (\resp. $\wh\cM^*$). We hence have $\wh{\cT_\cM}=\cT_{\wh\cM}$. 

Remark also that the dual $\wh\cD$-module of $\wh\cM^{\ir}$ is a $\oxdh[*D]$-module (\cite[chap\ptbl III, rem\ptbl (2.4)]{Malgrange91}). This implies that $\wh\cM_{\min}^{\ir}=\wh\cM^{\ir}$ and $\cT_{\wh\cM}=\cT_{\wh\cM^{\reg}}$. It follows that $\cM_{\min}$ is equal to the kernel of the composed morphism
$$
\cM\longrightarrow \wh\cM\longrightarrow \cT_{\wh\cM}=\cT_{\wh\cM^{\reg}}.
$$

We have recalled above that when $\cM$ is regular, we have $\DR(\cM_{\min})=j_*M$.

On the other hand, if $\cM$ is purely irregular (\ie. if $\wh\cM^{\reg}=0$), then $\cM_{\min}=\cM$ by the remark above, and consequently $\DR(\cM_{\min})=\DR(\cM)$. In this case, $\DR(\cM_{\min})$ is far from being concentrated in degree $0$ because, at each point of $D$, the difference $\dim\cH^1-\dim\cH^0$ is equal to the irregularity number of $\cM$ at this point.

\subsection{}
Let $(\Db_X^\bbullet,d)$ be the complex of currents on $X$. The natural morphism $(\Omega_X^\bbullet,d)\rightarrow (\Db_X^\bbullet,d)$ is a quasi-isomorphism. As $\cM_{(\min)}$ (\ie. $\cM$ or $\cM_{\min}$) is $\cO_X$-flat, the complex $(\Db_X^{k,\bbullet}\otimes_{\cO_X}\cM_{(\min)},\ov\partial)$ is a resolution of $\Omega_{X}^{k}\otimes _{\cO_X}\cM_{(\min)}$, and, taking the single complex associated with the double complex $(\Db_X^{\bbullet,\bbullet}\otimes_{\cO_X}\cM_{(\min)},\nabla,\ov\partial)$ we conclude that
$$
\DR(\cM_{(\min)})= (\Omega_X^\bbullet\otimes_{\cO_X}\cM_{(\min)},\nabla)\isom (\Db_X^\bbullet\otimes_{\cO_X}\cM_{(\min)},\nabla+\ov\partial).
$$

Fix a complete metric on $X^*$ as in \T\ref{sec:statement}. If a metric is given on the bundle $\cM_{|X^*}$, the sheaf $\cL_{(2)}^{i}(\cM)$ is the sheaf on $X$ which sections on an open set $U$ of $X$ are $i$-currents $\tau$ on $U^*$ with values in $\cM$, such that $\tau$ and $(\nabla+\ov\partial)\tau$ are in $L^2$ on each compact set of $U$, with respect to the chosen metrics. If the matrix of the metric on $\cM_{|X^{*}}$ in some (or any) local meromorphic basis of $\cM$ has moderate growth along $D$, the complex $\cL_{(2)}^{\bbullet}(\cM)$ is naturally a subcomplex of
$$
\Db_X^{\mod 0,\bbullet}\otimes \cM=\Db_X^\bbullet[*D]\otimes \cM=\Db_X^{\bbullet}\otimes \cM,
$$
where $\Db_X^{\mod 0}=\textrm{image}\lefcro\Db_X\rightarrow \Db_{X^*}\rigcro$.

Denote by $\cC$ the cone of $\cL_{(2)}^{\bbullet}(\cM)\to\Db_X^\bbullet\otimes \cM$. To prove theorem \ref{th:poincare} it is enough to construct a morphism $\cC\to\DR(\cT_{\cM})$ in $D^{b}(\CC_{X})$ such that the diagram
$$
\begin{array}{c}
\xymatrix{
\Db_X^\bbullet\otimes \cM\ar[r]&\cC\ar@{-->}[d]\\
\DR(\cM)\ar[u]_{\wr}\ar[r]&\DR(\cT_{\cM})
}
\end{array}
\leqno{(*)}
$$
commutes and to show that it is an isomorphism. This is now a local problem on $X$, so we may (and will) assume in the proof of theorem \ref{th:poincare} that $X=\Delta$ and $D=\{0\}$. It is also enough to prove the result for germs at the origin, because it is clear outside the origin.

\subsection{}\label{subsec:tilde}
Let $e:\wt\Delta=[0,r_0[{}\times S^1\rightarrow \Delta$ be the real blow-up of $\Delta$ at the origin (polar coordinates) and let $\cA_{\wt\Delta}$ be the subsheaf of $\cC_{\wt\Delta}^{\infty}$ defined as the kernel of the Cauchy-Riemann operator $\ov z\partial_{\ov z}$. Put $\wt\cM=\cA_{\wt\Delta}\otimes_{e^{-1}\cO_\Delta}e^{-1}\cM$. As $\partial_z$ operates on sections of $\cA_{\wt\Delta}$, one may extend the action of $\nabla$ to $\wt\cM$ in a natural way.

The de~Rham complex $\DR(\wt\cM)$ is the complex $\cA_{\wt\Delta}\otimes (e^{-1}\Omega_\Delta^\bbullet\otimes \cM)$ with its natural differential.

The complex of currents with moderate growth
\[
\Db_{\wt\Delta}^{\mod0,\bbullet}\otimes_{e^{-1}\cO_\Delta}e^{-1}\cM=\Db_{\wt\Delta}^{\mod0,\bbullet}\otimes_{\cA_{\wt\Delta}}\wt\cM
\]
is also defined, because $\partial_z$ and $\partial_{\ov z}$ operate on $\Db_{\wt\Delta}^{\mod0}=\textrm{image}\lefcro \Db_{\wt\Delta}\rightarrow \Db_{\Delta^*}\rigcro$.

The complex $\cL_{(2)}^{\bbullet}(\wt\cM)$ can be defined in a similar way on $\wt\Delta$, and for the metric constructed in \T\ref{sec:construction}, we have a natural inclusion $\cL_{(2)}^{\bbullet}(\wt\cM)\hookrightarrow \Db_{\wt\Delta}^{\mod0,\bbullet}\otimes_{\cA_{\wt\Delta}} \wt\cM$.

Since the terms of these complexes are $c$-soft, and since we evidently have $e_*\cL_{(2)}^{\bbullet}(\wt\cM)=\cL_{(2)}^{\bbullet}(\cM)$, we conclude that this natural inclusion becomes the natural inclusion $$\cL_{(2)}^{\bbullet}(\cM)\hookrightarrow \Db_{\Delta}^{\mod0,\bbullet}\otimes_{\cO_{\Delta}}\cM= \Db_{\Delta}^{\bbullet}\otimes_{\cO_{\Delta}}\cM,
$$
after applying the functor $\bR e_*$.

\subsection{}\label{subsec:lift}
Let $(\cM,\nabla)$ be a free $\cO_\Delta[z^{-1}]$-module of rank $d$ with connection and fix a formal model $(\cM^{\el},\nabla^{\el})$ as in lemma \ref{lem:modele} and an isomorphism $\wh\lambda:(\wh\cM^{\el},\wh\nabla^{\el})\rightarrow (\wh\cM,\wh\nabla)$. Let $\pi:\Delta_q\rightarrow \Delta$ be a cyclic covering such that $\pi^+(\cM^{\el},\nabla^{\el})$ is elementary. Then, for any $\theta^o\in S^1$, denoting by $\cA_{\theta^o}$ the germ of $\cA_{\wt\Delta}$ at $(0,\theta^o)\in[0,r_0[{}\times S^1$, there exists $\lambda_{\theta^o}$ such that the diagram
\[
\xymatrix{
(\cA_{\theta^o}\otimes \cM^{\el},\nabla^{\el})\ar[d]\ar[rr]^{\lambda_{\theta^o}}&&(\cA_{\theta^o}\otimes \cM,\nabla)\ar[d]\\
(\wh\cM^{\el},\wh\nabla^{\el})\ar[rr]^{\wh\lambda}&& (\wh\cM,\wh\nabla)
}
\]
commutes, where the vertical maps are induced by the Taylor expansion map $\cA_{\theta^o}\rightarrow \CC\lcr z\rcr$.

This is known (\cite{Wasow65,Sibuya90bis,Malgrange91}) when $(\wh\cM^{\el},\wh\nabla^{\el})$ is elementary. In general, remark that $\pi$ induces a covering map $\wt\pi:\wt\Delta_q\rightarrow \wt\Delta$, with $\wt\Delta_q=[0,r_{0}^{1/q}[{}\times S^1$, hence an isomorphism $\wt\pi_{\theta^o}^*:\cA_{\theta^o}\isom\cA_{\theta^o/q}$ for any possible choice of $\theta^o/q$, and one may define $\lambda_{\theta^o}$ as
\[
(\wt\pi_{\theta^o}^{-1}\otimes \id)^{-1}\circ\lambda_{\theta^o/q}\circ(\wt\pi_{\theta^o}^{-1}\otimes \id).
\]

It follows that there exists a covering of $S^1$ by open intervals $I$ such that all triple intersections are empty and all intersections are connected or empty, and isomorphisms
\[
\lambda_I:(\cA\otimes \cM^{\el},\nabla^{\el})_{|[0,\varepsilon[{}\times I}\isom (\cA\otimes \cM,\nabla)_{|[0,\varepsilon[{}\times I}
\]
lifting $\wh\lambda$, for $\varepsilon>0$ small enough. The Stokes automorphisms $\lambda_J\circ\lambda_{I}^{-1}=\id+\mu_{I,J}$ on $[0,\varepsilon[{}\times(I\cap J)$ are such that $\mu_{I,J}$ is nilpotent and its matrix in any meromorphic basis of $\cM$ is $C^\infty$-flat along $\{0\}\times (I\cap J)$; more precisely, the entries have an exponential decay along $\{0\}\times (I\cap J)$. Let us prove this assertion.

It is possible to write
\[
\pi^{-1}(I)=I_1\cup\cdots\cup I_q,\quad\pi^{-1}(J)=J_1\cup\cdots\cup J_q,
\]
a disjoint union of open intervals of $S^1$ such that $I_i\cap I_j\neq\emptyset$ for $i\neq j$ and $\wt\pi:I_k\cap J_k\rightarrow I\cap J$ is a diffeomorphism for any $k=1,\ldots ,q$. Put $\lambda_{I_k}=(\wt\pi_{|I_k}\otimes \id)^{-1}\circ\lambda_I\circ (\wt\pi_{|I_k}\otimes \id)$. It is known that the matrix of $\lambda_{I_k}\circ\lambda_{J_k}^{-1}$ in any meromorphic basis of $\pi^*\cM^{\el}$ is equal to $\id+N_{I_k,J_k}$ with $N_{I_k,J_k}$ nilpotent and exponentially decreasing at $\{0\}\times (I_k\cap J_k)$. Hence the same property holds for $\lambda_I\circ\lambda_{J}^{-1}$ on $\cA\otimes \cM^{\el}_{|[0,\varepsilon[{}\times (I\cap J)}$. This is also true for the matrix of $\lambda_I\circ\lambda_{J}^{-1}$ in any $\cA$-basis of $\cA\otimes \cM^{\el}_{|[0,\varepsilon[{}\times(I\cap J)}$, hence for the matrix of $\lambda_J\circ\lambda_{I}^{-1}$ in any $\cA$-basis of $\cA\otimes \cM_{|[0,\varepsilon[{}\times(I\cap J)}$.\hfill\qed

\subsection{}\label{subsec:determinant}
Remark that, for any rank-one meromorphic connection $(\cL,\nabla)$ on $\Delta$ with pole at $0$, the natural connection on $\cHom_{\cO_\Delta}(\cL,\cL)$ is regular at $0$. It follows that for any meromorphic connection equipped with a formal isomorphism $\wh\lambda$ as above, the isomorphism $\det\wh\lambda$ can be lifted to an isomorphism $\det(\cM^{\el},\nabla^{\el})\isom\det(\cM,\nabla)$. In particular, given any basis $\bme^{\el}$ of $\cM^{\el}$, there exists a basis $\bme$ of $\cM$ such that the matrix $\wh\Lambda$ of $\wh\lambda$ in these bases satisfies $\det\wh\Lambda=1$.

Let $\Lambda_I$ be the matrix of $\lambda_I$ in these bases. Then $\det\Lambda_I=\det\Lambda_J$ on $[0,\varepsilon[{}\times I\cap J$ as $\mu_{I,J}$ is nilpotent. Hence $\det\Lambda_I$ is a meromorphic function independent of $I$. Its asymptotic expansion at $0$ being equal to $1$, it is equal to $1$. Therefore, in the bases $\bme^{\el}$ and $\bme$ as above, we have $\det\Lambda_I=1$ for all $I$.

\section{Construction of the metric $k$}\label{sec:construction}
The metric $k$ will be constructed in the neighbourhood of each point of $D$ and then will be extended to $X^*$ in a $C^\infty$ way using a partition of unity.

The metric $k$ on $\cC^\infty_{\Delta^*}\otimes \cM_{|\Delta^*}$ will be first defined for elementary meromorphic connections, then  for any connection.

\subsection{}\label{subsec:treselem}
See also \cite[\T5]{Simpson90} and \cite[\T11]{Biquard97}.

Let $\alpha\in\CC$, $\varphi\in 1/z\cdot\CC[1/z]$ and put $\alpha'=\reel\alpha$, $\alpha''=\im\alpha$.

Let $Y$ be a nilpotent $d\times d$ matrix, $(Y,X,H)$ a $\sld$-triple and $\Crochet$ a positive definite hermitian form such that $Y^*=X$, $X^*=Y$ and $H^*=H$. Fix an orthonormal basis $\bme^o$ of $V=\CC^d$, made with eigenvectors of $H$ and denote by $w_j$ the eigenvalue of $H$ corresponding to $e_j^o$, with $w_j\in\ZZ$. We will denote by $Y,X,H$ the matrices of $Y,X,H$ in the basis $\bme^o$.

Put on the trivial bundle $\ccV=\cC^\infty_\Delta\otimes_\CC V$ on the disc $\Delta=\{z\mid\module{z}<r_0\}$ the connection $D$ such that, in the basis $\bme=1\otimes \bme^o$, one has
\begin{eqnarray*}
D^{0,1}\bme&=&0\\
D^{1,0}\bme&=&\bme\cdot(Y+(-\alpha+z\varphi'(z))\id)\otimes \frac{dz}{z}.
\end{eqnarray*}
The family $\wt{\bme}$ defined by
\begin{eqnarray*}
\wt e_j&=&e^{-\varphi(z)}z^\alpha\exp(-Y\log z)\cdot e_j\\
&=&e^{-\varphi(z)}z^\alpha\sum_{k\geq 0}(-\log z)^k\cdot Y^ke_j
\end{eqnarray*}
is a multivalued horizontal basis of $\ccV_{|\Delta^*}$.

\medskip
Put $a(z)=\module{\log z\ov z}$ and let $P(z)$ be the matrix
\begin{eqnarray*}
P(z)&=&\module{z}^{\alpha'}a(z)^{-H/2}e^X.
\end{eqnarray*}

Consider the basis $\varepsilong$ of $\ccV_{|\Delta^*}$ defined as
\begin{eqnarray*}
\varepsilong&=&\bme\cdot P(z).
\end{eqnarray*}

The metric $k$ will be the metric on $\ccV_{|\Delta^*}$ such that the basis $\varepsilong$ is orthonormal.

The metric $k_1$ for which an orthonormal basis is $\bme\cdot \module{z}^{\alpha'}a(z)^{-H/2}$ will also be useful. Of course, $k$ and $k_1$ are mutually bounded. Notice that $\bme$ is $k_1$-orthogonal and that
\begin{eqnarray*}
k_1(e_j,e_j)&=&\module{z}^{-2\alpha'}a(z)^{w_j};
\end{eqnarray*}
but the horizontal multivalued basis $\wt\bme$ is no more orthogonal; however, for any closed sector $\{z\in\Delta^*\mid \arg z\in[\theta_0,\theta_1]\}$ there exist two constants $C_1,C_2$ such that any branch of $\wt e_j$ satisfies on this sector
\[
C_1e^{-\reel\varphi(z)}a(z)^{w_j}\leq k_1(\wt e_j,\wt e_j)\leq C_2 e^{-\reel\varphi(z)}a(z)^{w_j}.
\]

\begin{remarque}
When $\varphi\equiv0$ and $\alpha$ is real, the metric $k_1$ is the metric used by Zucker in \cite[prop\ptbl 11.3]{Zucker79}. His argument also gives a proof of theorem \ref{th:poincare} for any $\alpha\in\CC$ (and $\varphi\equiv0$), when one uses the metric $k_1$ to compute the $L^2$ complex. Notice that, as $k_1$ and $k$ are mutually bounded, they define the same $L^2$ complex.
\end{remarque}

\begin{lemme}
The curvature $R_k$ of the metric $k$ is ``acceptable'' in the sense of \cite{Simpson90}; more precisely it satisfies $\norme{R_k}_k\leq \dfrac{C}{\module{z}^2\module{\log z\ov z}^2}$ with $C>0$, the norm being computed with the standard Euclidian metric on $\Delta$.
\end{lemme}

\begin{proof}
Recall that for $\ell\in\ZZ$ and $z\in\Delta$ we have
\[
z\partial_za(z)^\ell=\ov z\partial_{\ov z}a(z)^\ell=-\ell a(z)^{\ell-1}.
\]
We will use the following identities:
\begin{eqnarray*}
a^{\pm H/2}Ya^{\mp H/2}&=&a^{\mp1}Y\\
a^{\pm H/2}Xa^{\mp H/2}&=&a^{\pm1}X\\
e^YHe^{-Y}&=&H+2Y\\
e^XHe^{-X}&=&H-2X\\
e^XYe^{-X}&=& Y+H-X.
\end{eqnarray*}
The matrix $K$ of $k$ in the basis $\bme$ is
\[
K=(P^{-1})^*P^{-1}=\module{z}^{-2\alpha'}a^{H/2}e^{-Y}e^{-X}a^{H/2}
\]
and the matrix of the metric connection $D_k=D_k^{1,0}+\ov \partial$ is
\begin{eqnarray*}
M_k^{1,0}\frac{dz}{z}&=&\lefcro-\alpha'\id+K^{-1}(z\partial_z+\alpha')K\rigcro\frac{dz}{z}.
\end{eqnarray*}
The computation gives
\begin{eqnarray*}
M_k^{1,0}&=&-\alpha'\id-Y-\frac{2H}{a}+\frac{2X}{a^2}.
\end{eqnarray*}
In the basis $\bme$, the curvature $R_k$ has matrix $R_k^{1,1}\,\dfrac{dz\wedge d\ov z}{\module{z}^2}$ with
\[
R_k^{1,1}=-\ov z\partial_{\ov z}M_k^{1,0}=\frac{2H}{a^2}-\frac{4X}{a^3}.
\]
The matrix of $R_k$ in the orthonormal basis $\varepsilong$ is thus
\begin{eqnarray*}
P^{-1}R_k^{1,1}P\,\frac{dz\wedge d\ov z}{\module{z}^2}&=& \frac{2H\,dz\wedge d\ov z}{\module{z}^2a^2}
\end{eqnarray*}
and the ``acceptability'' is clear.
\end{proof}

\begin{lemme}
The pseudo-curvature $G_k$ of $(\ccV_{|\Delta^*},D,k)$ is zero.
\end{lemme}

\begin{proof}
Denote by $M^{1,0}dz/z$ (\resp.\ $M^{0,1}d\ov z/\ov z$) the matrix of $D^{1,0}$ (\resp.\ $D^{0,1}$) in the basis $\varepsilong$. We have
\begin{eqnarray*}
M^{1,0}&=&(-\alpha+z\varphi'(z))\id+P^{-1}YP+P^{-1}z\partial _zP,\\
M^{0,1}&=&P^{-1}\ov z\partial_{\ov z}P.
\end{eqnarray*}
Put
\[
\Theta=\frac12(M^{1,0}+M^{0,1*}),\quad
\Theta^*=\frac12(M^{1,0*}+M^{0,1})
\]
and $\theta=\Theta dz/z$, $\theta^*=\Theta^*d\ov z/\ov z$. The matrix $N^{0,1}d\ov z/\ov z$ of the operator $\ov\partial _E\defin D^{0,1}-\theta^*$ in the basis $\varepsilong$ is $N^{0,1}=M^{0,1}-\Theta^*$. The pseudo-curvature $G_k\defin\ov\partial _E(\theta)$ has matrix
\[
\lefpar \ov z\partial _{\ov z}(\Theta)+[N^{0,1},\Theta]\rigpar \frac{d\ov z}{\ov z}\wedge\frac{dz}{z}.
\]
Put $\Theta=\Theta_{\Reg}+\dfrac{z\varphi'(z)}{2}\id$ and $N_{\Reg}^{0,1}=M^{0,1}-\Theta_{\Reg}^{*}$. Then
\begin{eqnarray*}
\ov z\partial _{\ov z}(\Theta)+[N^{0,1},\Theta]&=&\ov z\partial _{\ov z}(\Theta_{\Reg})+[N_{\Reg}^{0,1},\Theta_{\Reg}].
\end{eqnarray*}
It is thus enough to prove the lemma when $\varphi\equiv0$. One shows that in the basis $\bme\cdot Q(z)$ with $Q(z)=\module{z}^{i\alpha''}a^{-H/2}e^Xa^{H/2}$, the matrix of $\ov\partial_{E\Reg}$ is zero and the matrix of $\theta_{\Reg}$ is equal to
\[
\lefpar \dfrac{-i\alpha''}{2}\id+Y\rigpar \frac{dz}{z}
\]
which is holomorphic, hence $\ov\partial_{E\Reg}(\theta_{\Reg})=0$.
\end{proof}

\begin{remarque}
The notation $\ov\partial_E$ is taken from \cite{Simpson90}. The holomorphic bundle $E$ is the bundle generated by the basis $\bme\cdot Q(z)e_{}^{(\ov \varphi-\varphi)/2}$. The computation above associates to $(\ccV,D)$ and the metric $k$ a \emph{Higgs bundle} $(E,\theta)$.
\end{remarque}

\subsection{}\label{subsec:elem}
Let now $V$ be a $\CC$-vector space of dimension $d$, equipped with an automorphism $T$ and with a $T$-stable decomposition $V=\oplus_{\varphi\in\Phi}V_\varphi$ indexed by a finite subset $\Phi\subset z^{-1}\CC[z^{-1}]$. Let $T=T_sT_u$ be the decomposition of $T$ into its semi-simple and unipotent part and put $Y=-\dfrac{1}{2i\pi}\log T_u$. Fix an $\sld$-triple $(Y,X,H)$ compatible with $T_s$ and with the $\Phi$-decomposition of $V$, and fix a positive definite hermitian form $\Crochet$ such that
\begin{enumerate}
\item[\refstepcounter{theoreme}\label{221}(\thetheoreme)]
The double decomposition $V=\oplus_{\varphi,\alpha}V_{\varphi,\alpha}$ with respect to $\Phi$ and to $\dfrac{1}{2i\pi}\log$ of the eigenvalues of $T$ is orthogonal,
\item[\refstepcounter{theoreme}\label{222}(\thetheoreme)]
$Y^*=X$, $X^*=Y$ and $H^*=H$.
\end{enumerate}
Fix also an orthonormal basis $\bme^o$ of eigenvectors of $H$ adapted to this double decomposition.

Let $\bme=1\otimes \bme^o$ be the corresponding basis of $\ccV_\cO=\cO_\Delta\otimes V$. Put on the meromorphic bundle $\ccV_\cO[z^{-1}]$ the connection $\nabla$ which is the direct sum of those used in \T\ref{subsec:treselem}. This meromorphic connection is elementary and any elementary connection has this form. 

Define the metric $k$ in such a way that $\varepsilong=\bme\cdot P(z)$ is a $k$-orthonormal $\ccV_{|\Delta^*}$-basis, where $P$ is a block-diagonal matrix, the blocks being those of \T\ref{subsec:treselem}. Define the metric $k_1$ in a similar way. The matrix $P$ can be written as $P(z)=\delta(z)a(z)^{-H/2}e^X$, where $\delta(z)$ is a diagonal matrix with diagonal entries of the form $\module{z}^{\alpha'}$ and commuting with $Y,X,H$.

Clearly, the metric $k$ satisfies both lemmas of \T\ref{subsec:treselem}.

\subsection{}\label{subsec:elemramif}
Keep the same notation as above. Let $q\in\NN^*$ and put $\zeta=e^{2i\pi/q}$. Let $(\cM,\nabla)$ be a meromorphic bundle with connection on $\Delta$, with pole at $0$. Let $\pi:t\mapsto t^q=z$ be a cyclic ramified covering from a disc $\Delta_q$ to $\Delta$ and assume that $\pi^+(\cM,\nabla)$ is elementary, \ie. can be described as in \T\ref{subsec:elem}. As $\pi=\pi\circ\sigma$, we moreover have an isomorphism $\lambda_\sigma:\pi^+(\cM,\nabla)\isom\sigma^+(\pi^+(\cM,\nabla))$. In terms of the data of \T\ref{subsec:elem}, we are given an automorphism of finite order $q$:
\[
\lambda_{\sigma}^{o}:V\longrightarrow V
\]
which commutes with $T_u$ and is $\sigma$-compatible with the double decomposition of $V$, \ie., $\lambda_{\sigma}^{o}(V_{\varphi,\alpha})=V_{\varphi\circ\sigma,\alpha+1/q}$, where $\alpha+1/q$ is taken modulo $1$. We will choose $H$ and $X$ so that $\lambda_{\sigma}^{o}$ also commutes with the $\sld$-triple. It defines an isomorphism
\begin{eqnarray*}
\ccV&\MRE{\lambda_\sigma}&\sigma^*\ccV\\
f(t)\otimes v^o&\DMRE{}&f(\zeta^{-1}t)\otimes \lambda_{\sigma}^{o}(v^o).
\end{eqnarray*}
Fix a positive definite hermitian form $\Crochet$ on $V$ which satisfies \eqref{221}, \eqref{222} and
\begin{enumerate}
\item[\refstepcounter{theoreme}\label{231}(\thetheoreme)]
\quad$\Crochet$ is $\lambda_{\sigma}^{o}$-invariant.
\end{enumerate}
Let $\bme^o$ be an orthonormal basis of $V$ made with eigenvectors of $H$ and compatible with the double decomposition, and put $\bme=1\otimes \bme^o$, $\varepsilong=\bme\cdot P(t)$ with $P$ as in \T\ref{subsec:elem}. We thus get a metric $k_q$ on $\ccV_{|\Delta_q^*}$.

\begin{lemme}
The metric $k_q$ is compatible with $\lambda_\sigma$, \ie.,
\[
k_{q,\sigma}(v,v')\defin\sigma^*k_q(\lambda_\sigma(v),\lambda_\sigma(v'))= k_q(v,v')
\]
for local sections $v,v'$ of $\ccV$.
\end{lemme}

\begin{proof}
Working in the basis $\bme$ and denoting by $\Lambda_{\sigma}^{o}$ the matrix of $\lambda_{\sigma}^{o}$ in this basis, we see that the matrix of the hermitian form $\sigma^*k_q(\lambda_\sigma(v),\lambda_\sigma(v'))$ is equal to $\Lambda_{\sigma}^{o*}\cdot K_q\circ\sigma\cdot\Lambda_{\sigma}^{o}$, if $K_q=(P^{-1})^*P^{-1}$ denotes the matrix of $k_q$. As $\Lambda_{\sigma}^{o}$ is orthogonal and commutes with $H$ and $X$ and as $P\circ\sigma=P$, this matrix is equal to $K_q$.
\end{proof}

We deduce from this lemma that there exists a metric $k$ on $\cM_{|\Delta^*}$ such that the metric $k_q$ on $\ccV_{|\Delta_q^*}$ is pulled-back from $k$ by $\pi$. In particular, the curvature and the pseudo-curvature are also pulled-back, and both lemmas of \T\ref{subsec:treselem} are satisfied by $k$ on $\cM_{|\Delta^*}$ as soon as they are satisfied by $k_q$ on $\pi^*\ccV_{|\Delta_q^*}$.

\begin{Remarques}
\begin{enumerate}
\item
The determinant $\det K_q$ is equal to $(\det P)^{-2}= \module{t}^{2\beta'}$ for some $\beta'\in\RR$, as $\tr H=0$ and $X$ is nilpotent.
\item
The metric $k_{1,q}$ is also compatible with $\lambda_\sigma$, hence defines a metric $k_1$ on $\cM_{|\Delta^*}$. Moreover, $k$ and $k_1$ are mutually bounded.
\end{enumerate}
\end{Remarques}

\subsection{}\label{subsec:decform}
Let $(\cM,\nabla)$ be a meromorphic bundle with connection on $\Delta$ having a singularity at $0$ only. Let $(\cM^{\el},\nabla^{\el})$ be a formal model for it, in the sense of lemma \ref{lem:modele}. Fix a formal isomorphism $\wh\lambda$ and liftings $\lambda_I$ as in \T\ref{subsec:lift}.

Index the family of intervals by $\ZZ/N\ZZ$, so that $I_\ell\cap I_{\ell+1}\neq\emptyset$ and $I_\ell\cap I_{\ell+k}=\emptyset$ if $k\neq-1,0,1$ for any $\ell$. Let $\chi_\ell\in\cC^\infty(I_\ell\cap I_{\ell+1})$ be such that $\chi_\ell\equiv0$ on the $\ell$-side and $\chi_\ell\equiv1$ on the $\ell+1$ side of $I_\ell\cap I_{\ell+1}$.

We may now define an isomorphism
\[
\wt\lambda:\cC_{[0,\varepsilon[{}\times S^1}^{\infty}\otimes \cM'\isom\cC_{[0,\varepsilon[{}\times S^1}^{\infty}\otimes \cM
\]
by the formula
\begin{eqnarray*}
\wt\lambda&=&
\begin{cases}
\lambda_{I_\ell}\text{ (also denoted $\lambda_\ell$)}&\text{on } [0,\varepsilon[{}\times (I_\ell-\ccup_{\ell'\neq\ell}I_{\ell'}),\\
(\id+\chi_\ell(\theta)\mu_{\ell,\ell+1})\circ\lambda_\ell&\text{on } [0,\varepsilon[{}\times (I_\ell\cap I_{\ell+1}).
\end{cases}
\end{eqnarray*}

Denote by $k^{\el}$ the metric on $\cM_{|\Delta^*}^{\el}$ constructed in \T\ref{subsec:decform} and let $k$ be the metric on $\cM_{|\Delta^*}$ obtained by pushing $k^{\el}$ by $\wt\lambda$. Define $k_1$ in a similar way, pushing $k_1^{\el}$ by $\wt\lambda$.

Let us show that $k$ satisfies (1.1), (1.2) and (1.3). Assume first, for simplicity, that $q=1$.

Denote now by $\bme^{\el}$ the basis of $\cM^{\el}$ used in \T\ref{subsec:elem}, and by $\varepsilong^{\el}$ the corresponding orthonormal basis; we have $\varepsilong^{\el}=\bme^{\el}\cdot P(z)$ where $P$ is a block-diagonal matrix as in \T\ref{subsec:elem}.

Put $\wt\bme=\wt\lambda(\bme^{\el})$. Then $\varepsilong=\wt\lambda(\varepsilong^{\el})=\wt\bme\cdot\wt P$ is an orthonormal basis for $k$.

In a neighbourhood of the sets $]0,\varepsilon[{}\times (I_\ell-\cup_{\ell'\neq\ell}I_{\ell'})$, the basis $\wt\bme$ is holomorphic and we may compute the curvature and the pseudo-curvature of $k$ exactly as in \T\ref{subsec:treselem}.

On $]0,\varepsilon[{}\times (I_\ell\cap I_{\ell+1})$, denote by $\bme^{(\ell)}$ the holomorphic basis $\lambda_\ell(\bme^{\el})$. We then have
\begin{eqnarray*}
\varepsilong&=&\bme^{(\ell)}(\id+M_{\ell,\ell+1})\cdot P,
\end{eqnarray*}
where $M_{\ell,\ell+1}$ is a nilpotent matrix with exponentially decreasing $C^\infty$ coefficients. The curvature and the pseudo-curvature of $k$ can be expressed with the same formulas as in \T\ref{subsec:treselem}, adding a perturbation term which is asymptotically equal to $0$ along $\{0\}\times S^1$.

Therefore, properties (1.1), (1.2) and (1.3) are satisfied for $(\cM_{|\Delta^*},\nabla,k)$. For $q\geq 2$, argue similarly with $\pi^+\cM^{\el}$ and $\pi^+\cM$.

\begin{remarque}\label{rem:last}
If $\bme$ is any meromorphic basis of $\cM$, the norm of $\det\bme$ for the metric on $\det\cM_{|\Delta^*}$ induced by $k$ or $k_1$ is equal to $\module{z}^{2\beta'}$ for some $\beta'\in\RR$ and some choice of coordinate on $\Delta$: indeed, this is true if $\bme$ is chosen so that the matrix of $\Lambda_\ell$ in the bases $\bme^{\el}$ and $\bme$ has determinant equal to $1$ for all $\ell$ (see \T\ref{subsec:determinant}); hence this is true for any meromorphic basis after a suitable change of coordinate.
\end{remarque}

\subsection{}\label{subsec:degre}
Let us now come back to the global setting. In order to show that the metric $k$ satisfies (1.4), we will need the lemma below. Let $\cL$ be a rank-one meromorphic bundle on $X$, with poles on $D$, equipped with a connection $\nabla:\cL\rightarrow \cL\otimes_{\cO_X}\Omega_X^1$. Let $x^o\in D$ and let $e$ be a local section of $\cL_{x^o}$ for which, if $z$ is a local coordinate at $x^o$, one has $\nabla e=e\otimes \omega$, with $\omega=(\alpha_{x^o}+z\varphi'_{x^o}(z))dz/z$, $\varphi_{x^o}\in z^{-1}\CC[z^{-1}]$ and $\alpha_{x^o}\in\CC$. Let $k$ be a metric on $\cL_{|X^*}$ which satisfies, near any $x^o\in D$,
\begin{eqnarray*}
k(e,e)&=&\module{z}^{2\reel(\alpha_{x^o})}.
\end{eqnarray*}

\begin{lemme}
With these assumptions, we have $\deg(\cL_{|X^*},k)=0$.
\end{lemme}

\begin{proof}
Let $L$ be the rank-one $\cO_X$-submodule of $\cL$ which is equal to $\cL_{|X^*}$ on $X^*$ and which is generated by the section $e$ near $x^o$, for any $x^o\in D$. The residue formula and the fact that $\deg L$ is an integer give
\[
\deg L=\sum_{x\in D}\res_x\nabla=\sum_{x\in D}\alpha_x=\sum_{x\in D}\reel(\alpha_x).
\]
Indeed, if $L$ is trivial, this is the usual residue formula for meromorphic differential forms; if $L$ is not trivial, there exists on its dual $L^*$ a logarithmic connection which satisfies the residue formula; the residue formulas for the trivial bundle $L\otimes L^*$ and for $L^*$ give the formula for $L$.

On the other hand, in a punctured neighbourhood of $x^o\in D$, the $(1,1)$-form $\ov\partial \partial \log k(e,e)$ is identically $0$. The curvature of $k$ exists as a $(1,1)$-current and the degree of $L$ can be computed with this current. Taking into account the Dirac currents at each $x\in D$, one gets
\begin{eqnarray*}
\deg L&=&\sum_{x\in D}\reel(\alpha_x)+\deg(\cL_{|X^*},k).
\end{eqnarray*}
The lemma follows.
\end{proof}

According to remark \ref{rem:last}, the metric constructed in \T\ref{subsec:decform} on $\det\cM_{|X^*}$ satisfies the assumptions of the previous lemma. This shows that $\deg(\cM_{|X^*},k)=0$.

\section{Proof of theorem \ref{th:poincare} (first part)}\label{sec:firstpart}
We will use polar coordinates $(r,\theta)$ on $[0,r_0[{}\times S^1$, with $r_0<1$. We put on $]0,r_0[{}\times S^1=\Delta^*$ the Poincar{\'e} metric and we denote by $d\vol$ its volume form. We have
\begin{enumerate}\setcounter{enumi}{-1}
\item
$f\in L^2(d\vol)\ssi \module{\log r}^{-1}f\in L^2(d\theta\,dr/r)$;
\item
$\omega=f\,dr/r+g\,d\theta\in L^2(d\vol)\ssi f\text{ and }g\in L^2(d\theta\,dr/r)$;
\item
$\eta=h\,d\theta\,dr/r\in L^2(d\vol)\ssi \module{\log r}h\in L^2(d\theta\,dr/r)$.
\end{enumerate}

Let $\varphi(z)=-\dfrac{a_\ell}{z^\ell}(1+z\psi(z))$ with $\ell\in\NN^*$, $a_\ell\in\CC$ and $\psi$ holomorphic on $\Delta$. Let $\beta\in\RR$ and $k\in\ZZ$ (for the application we will take $\beta=0$).

Let $L$ be a flat bundle of rank one on $\Delta^*$ equipped with a metric such that a multivalued horizontal section $\wt e$ has norm $\norme{\wt e}^2\sim r^{2\beta}\module{\log r}^k e^{-2\reel\varphi(z)}$. Let $\cL_{(2)}^\bbullet(L)$ be the $L^2$ complex on $[0,r_0[{}\times S^1$: a local section of $\cL_{(2)}^i(L)$ on an open set $\Omega$ of $[0,r_0[{}\times S^1$ can be written $\omega\otimes \wt e$, where $\omega$ is a $i$-current on $\Omega\cap\Delta^*$ such that the currents
$$
r^\beta\module{\log r}^{k/2}e^{-\reel\varphi(z)}\omega\qqbox{and} r^\beta\module{\log r}^{k/2}e^{-\reel\varphi(z)}d\omega
$$
are in $L^2(d\vol)$ on any compact set of $\Omega$.

\begin{lemme}\label{lem:rg1}
In this situation we have
\begin{enumerate}
\item
$\cH^i(\cL_{(2)}^\bbullet(L))=0$ for $i\geq 2$;
\item
$\cH^1(\cL_{(2)}^\bbullet(L))=0$ if $a_\ell\neq0$ or $\beta\neq0$.
\end{enumerate}
\end{lemme}

\begin{remarque*}
In the proof we will assume that $a_\ell\neq0$. If $a_\ell=0$ (and $\beta\neq0$ in (2)), the proof is analogous and even simpler. When $a_\ell\neq0$ we put $-a_\ell=\module{a_\ell}e^{i\tau}$. We then have
\begin{eqnarray*}
-\reel(\varphi)&=&\dfrac{\module{a_\ell}}{r^\ell}\lefpar \cos(\ell\theta-\tau)+r\delta_\varphi(r,\theta)\rigpar
\end{eqnarray*}
such that $\delta_\varphi$ is real analytic on $[0,r_0[{}\times S^1$; hence there exists $r_\varphi\in{}]0,r_0]$ such that, for any $r\in{}]0,r_\varphi[$, the sign of $\partial (-\reel\varphi)/\partial \theta$ is equal to the sign of $\sin(\tau-\ell\theta)$ and the sign of $\partial (-\reel\varphi)/\partial r$ is equal to the sign of $-\cos(\tau-\ell\theta)$.

The property that we will use in an essential way is that the map
\begin{eqnarray*}
(r,\theta)&\longmapsto&e^{-\reel\varphi(re^{i\theta})}
\end{eqnarray*}
has no critical point on $]0,r_0[{}\times S^1$, \ie. $\partial (-\reel\varphi)/\partial r$ and $\partial (-\reel\varphi)/\partial \theta$ do not both vanish as soon as $r_0$ is sufficiently small.
\end{remarque*}

\begin{proof}
Let $\Omega={}]0,r_1[{}\times {}]\theta_0,\theta_1[$ with $r_1\leq r_\varphi$. Let $\Omega''={}]0,r''_1[{}\times {}]\theta''_0,\theta''_1[$ containing $\ov\Omega\cap\Delta^*$. We assume that $]\theta''_0,\theta''_1[$ contains at most one zero of $\sin(\tau-\ell\theta)\cos(\tau-\ell\theta)$ and, if it contains one, this zero belongs to $]\theta_0,\theta_1[$.

Let $\omega=fdr+gd\theta$ (\resp. $\eta=h\,d\theta\,dr$) be a one-form (\resp. two-form) on $\Omega''$ such that $\omega\otimes \wt e\in L^2(\Omega'';L;d\vol)$ and $d\omega=0$ (\resp\ldots). We will show that there exists a function $u$ (\resp. a one-form $\xi$) on $\Omega$ such that $u\otimes \wt e\in L^2(\Omega;L;d\vol)$ and $du=\omega_{|\Omega}$ (\resp\ldots). This will give the vanishing of $\cH^1$ (\resp. $\cH^2$).

By assumption, $e^{-\reel\varphi}$ is monotonic with respect to $r$ or to $\theta$ on $\Omega''$.

Choose $\Omega'={}]0,r'_1[{}\times {}]\theta'_0,\theta'_1[$ with $\Omega\Subset\Omega'\Subset\Omega''$. There exists then (see \eg. \cite[prop\ptbl 12.2]{Demailly96}) a sequence $\varepsilon_n>0$ converging to $0$ and a sequence of one-forms $\omega_n$ with coefficients in $\cC_{0}^{\infty}(\Omega'')$ (\resp. of two-forms $\eta_n=h_nd\theta\,dr$) with support in $[\varepsilon_n,r''_1[{}\times {}]\theta''_0,\theta''_1[$ such that
$$
\omega_n\otimes \wt e \rightarrow \omega\otimes \wt e \;(\text{\resp. }\eta_n\otimes \wt e\rightarrow \eta\otimes \wt e)\qqbox{and}\quad d\omega_n\otimes \wt e \rightarrow 0\qqbox{in} L^2(\Omega';L;d\vol).
$$

\subsection*{First case: $e^{-\reel\varphi}$ is monotonic with respect to $\theta$ on $\Omega''$}
We assume here that $\sin(\tau-\ell\theta)$ does not vanish on $]\theta''_0,\theta''_1[$.

Let us begin with $\cH^1$. Let $\chi=\chi(\theta)\in\cC^\infty(S^1)$, with $\chi\equiv1$ on $[\theta_0,\theta_1]$ and $\chi\equiv0$ outside $]\theta'_0,\theta'_1[$. Put $\omega_n=f_ndr+g_nd\theta$ and define
\begin{eqnarray*}
u_n(r,\theta)&=&\begin{cases}
\dpl\int_{\theta'_0}^{\theta}\chi(t)g_n(r,t)dt&\text{if $e^{-\reel\varphi}$ is decreasing},\\
\dpl -\int_{\theta}^{\theta'_1}\chi(t)g_n(r,t)dt&\text{if $e^{-\reel\varphi}$ is increasing}.
\end{cases}
\end{eqnarray*}
We will consider the decreasing case for instance. One has a Hardy-type inequality (see for instance \cite[th\ptbl1.14]{O-K90})
\begin{multline*}
U_n(r)\defin\int_{\theta'_0}^{\theta'_1}\module{u_n(r,\theta)}^2 e^{-2\reel\varphi(re^{i\theta})}\,d\theta\\
\leq C_n(r)\int_{\theta'_0}^{\theta'_1}\module{\chi(\theta)g_n(r,\theta)}^2 e^{-2\reel\varphi(re^{i\theta})}\,d\theta\defin C_n(r)G_n(r)
\end{multline*}
with
\begin{eqnarray*}
C_n(r)&=&4\sup_{\theta\in[\theta'_0,\theta'_1]}\int_{\theta}^{\theta'_1}e^{-2\reel\varphi(re^{it})}\,dt\cdot \int_{\theta'_0}^{\theta}e^{2\reel\varphi(re^{it})}\,dt\\
&\leq &4\sup_\theta e^{-2\reel\varphi(re^{i\theta})}e^{2\reel\varphi(re^{i\theta})}(\theta'_{1}-\theta)(\theta-\theta'_{0})\\
&=&(\theta'_{1}-\theta'_{0})^{2}
\end{eqnarray*}
because $e^{-\reel\varphi}$ is decreasing on $[\theta'_0,\theta'_1]$. We hence have
\begin{multline*}
\norme{u_n\otimes \wt e}_{L^2(\Omega';L;d\vol)}^{2}= \int_{0}^{r_0}U_n(r)r^{2\beta-1}\module{\log r}^{k-2}\,dr\\
\leq (\theta'_{1}-\theta'_{0})^{2}\int_{0}^{r_0}G_n(r)r^{2\beta-1}\module{\log r}^{k}\,dr=4\norme{\chi g_nd\theta\otimes \wt e}_{L^2(\Omega';L;d\vol)}^{2}.
\end{multline*}
We thus deduce the existence of $u$ on $\Omega'$ such that
$$
r^{\beta-1/2}\module{\log r}^{k/2-1}e^{-\reel\varphi}\cdot u\in L^2(\Omega';d\theta\,dr)\qqbox{and}\frac{\partial u}{\partial \theta}=\chi g\text{ on } \Omega'.
$$

\begin{remarque*}
At this step, we can repeat the proof for $\cH^2$ and get the vanishing of $\cH^2$ except maybe at points where $\sin(\tau-\ell\theta)$ vanishes. Up to now we did not use the assumption $a_\ell\neq0$ or $\beta\neq0$. As indicated above, we will now assume that $a_\ell\neq0$, and leave the regular case $a_\ell=0$ to the reader.
\end{remarque*}

On the other hand we have
\begin{eqnarray*}
\frac{\partial u_n}{\partial r}\; =\; \int_{\theta'_0}^{\theta}\chi\frac{\partial g_n}{\partial r}\,dt &=& \chi(\theta)f_n(r,\theta)-\int_{\theta'_0}^{\theta}\chi'(t)f_n(r,t)\,dt - \int_{\theta'_0}^{\theta}\chi\lefpar \frac{\partial f_n}{\partial t}-\frac{\partial g_n}{\partial r}\rigpar dt.
\end{eqnarray*}

We deduce as above that the sequence $r^{\beta+1/2}\module{\log r}^{k/2}e^{-\reel\varphi}\cdot\partial u_n/\partial r$ has a limit in $L^2(\Omega';d\theta\,dr)$. Therefore, we have $u\otimes \wt e\in L^2(\Omega;L;d\vol)$ and $\omega-du=\ell(r)dr$ since $d\omega=0$.

We are hence reduced to the case where $\omega=f(r)dr$ with $\omega\otimes \wt e\in L^2(\Omega;L;d\vol)$.

Put $\psi(r)=\dpl\int_{\theta_0}^{\theta_1}e^{-2\reel\varphi(re^{i\theta})}d\theta$.

\begin{lemme}\label{lem:utile}
Assume that $\cos(\ell\theta-\tau)$ does not vanish on $]\theta''_0,\theta''_1[$. Then there exists, for all $N\in\ZZ$, a number $r_N(\varphi,\theta_0,\theta_1)>0$ such that the function $r^N\psi(r)$ is monotonic, as $\psi$ is, on $]0,r_N(\varphi,\theta_0,\theta_1)[$.
\end{lemme}

Let us keep this lemma for granted and let us end the proof of lemma \ref{lem:rg1} if $\cos(\ell\theta-\tau)$ and $\sin(\ell\theta-\tau)$ do not vanish on $[\theta_0,\theta_1]$. Then $r^{2\beta}\module{\log r}^k\psi(r)$ and $r^{2\beta}\module{\log r}^{k-2}\psi(r)$ are monotonic, as $\psi$ is, on $]0,r_N[$ for a suitable $N$. We may assume in the following that $r_1<r_N(\varphi,\theta_0,\theta_1)$ by choosing $\Omega$ smaller. Define
$$
u(r)=
\begin{cases}
\dpl\int_{0}^{r}f(\rho)\,d\rho&\text{if $\psi$ is decreasing,}\\
\dpl-\int_{r}^{r_1}f(\rho)\,d\rho&\text{if $\psi$ is increasing.}
\end{cases}
$$
We then have, using a Hardy-type inequality as above,
\begin{multline*}
\norme{u\otimes \wt e}^2=\int_{0}^{r_1}\module{u(r)}^2r^{2\beta-1}\module{\log r}^{k-2}\psi(r)\,dr\\
\leq C\int_{0}^{r_1}\module{f(r)}^2r^{2\beta+1}\module{\log r}^{k}\psi(r)\,dr=\norme{\omega\otimes \wt e}^2
\end{multline*}
with
\begin{eqnarray*}
C&=&4\sup_{\rho\in[0,r_1[}\int_{\rho}^{r_1}r^{2\beta-1}\module{\log r}^{k-2}\psi(r)\,dr\cdot \int_{0}^{\rho}\lefpar r^{2\beta+1}\module{\log r}^{k}\psi(r)\rigpar^{-1} dr,
\end{eqnarray*}
if for instance $\psi$ is decreasing. But we then have
\begin{eqnarray*}
\int_{\rho}^{r_1}r^{2\beta}\module{\log r}^{k-2}\psi(r)\,\frac{dr}{r}&\leq &\rho^{2\beta}\module{\log \rho}^{k-2}\psi(\rho),\\
\int_{0}^{\rho}\lefpar r^{2\beta}\module{\log r}^{k}\psi(r)\rigpar^{-1} \frac{dr}{r}&\leq &\lefpar \rho^{2\beta}\module{\log \rho}^{k}\psi(\rho)\rigpar ^{-1},
\end{eqnarray*}
hence $C\leq 4\sup_{\rho\in[0,r_1[}\rho(r_{1}-\rho)\module{\log \rho}^{-2}<+\infty$.

The case where $\psi$ is increasing is analogous.

\medskip
Assume now that $\cos(\ell\theta-\tau)$ vanishes at $\theta_{1/2}\in[\theta_0,\theta_1]$; from the assumption made at the beginning of the proof, we may even suppose that $\theta_{1/2}\in{}]\theta_0,\theta_1[$.

Consider the case where $\cos(\ell\theta-\tau)>0$ on $]\theta_0,\theta_{1/2}[$ (the other case is treated in a similar way). We then have $\sin(\ell\theta-\tau)\sim1$ near $\theta_{1/2}$, so $e^{-\reel\varphi}$ is decreasing with respect to $\theta$ on $]\theta_0,\theta_1[$. Let $\varepsilon>0$ be small enough and put
\begin{eqnarray*}
\psi_+(r)&=&\int_{\theta_0}^{\theta_{1/2}-\varepsilon}e^{-\reel\varphi}d\theta.
\end{eqnarray*}
It is decreasing with respect to $r$. Moreover, there exists $\kappa(\varepsilon)>0$ such that
\[
\psi_+(r)\leq \psi(r)\leq \kappa(\varepsilon) \psi_+(r).
\]
It will be clear from the proof of lemma \ref{lem:utile} that we may apply it to $\psi_+(r)$ (using there $\theta_{1/2}-\varepsilon$ instead of $\theta_1$). We may conclude the proof in this case, defining $u(r)$ by the formula
\begin{eqnarray*}
u(r)&=&\int_{0}^{r}f(\rho)d\rho.
\end{eqnarray*}
The constant $C$ is now bounded by $4\kappa(\varepsilon) \sup_{\rho\in[0,r_1[}\rho(r_{1}-\rho)\module{\log \rho}^{-2}<+\infty$.

\begin{proof}[Proof of lemma \ref{lem:utile}]
We have
\begin{eqnarray*}
N+r\frac{\partial (-\reel\varphi)}{\partial r}&=&\lefpar N-2\ell\frac{\module{a_\ell}}{r^\ell}\cos(\ell\theta-\tau)(1+r\cdot\mu(r,\theta))\rigpar
\end{eqnarray*}
with $\mu$ real analytic on $[0,r_1[{}\times S^1$ depending only on $\varphi$. Since $\module{\cos(\ell\theta-\tau)}$ is bounded from below on $[\theta_0,\theta_1]$, there exists $r_N(\varphi,\theta_0,\theta_1)$ such that, for $r<r_N(\varphi,\theta_0,\theta_1)$, the sign of
\begin{eqnarray*}
r\frac{\partial }{\partial r}(r^N\psi(r))&=&r^N\int_{\theta_0}^{\theta_1}\lefcro N-2\ell\frac{\module{a_\ell}}{r^\ell}\cos(\ell\theta-\tau)(1+r\cdot\mu(r,\theta)) \rigcro e^{-2\reel\varphi}\,d\theta
\end{eqnarray*}
is equal to the sign of $-\cos(\ell\theta-\tau)$.
\end{proof}

\subsection*{Second case: $e^{-\reel\varphi}$ is monotonic with respect to $r$ on $\Omega''$}
We assume now that $\cos(\ell\theta-\tau)$ does not vanish on $\Omega'$ but $\sin(\ell\theta-\tau)$ vanishes at one point of $]\theta_0,\theta_1[$.

Choose $\chi\in\cC_{0}^{\infty}([0,r'_1[)$ with $\chi\equiv1$ on $[0,r_1[$. Consider the same sequence $\omega_n$ as above and put
$$
h_n(r,\theta)=
\begin{cases}
\dpl\int_{0}^{r}f_n(\rho,\theta)\,d\rho&\text{if $e^{-\reel\varphi}$ is decreasing},\\
\dpl-\int_{r}^{r'_1}\chi(\rho)f_n(\rho,\theta)\,d\rho&\text{if $e^{-\reel\varphi}$ is increasing}.
\end{cases}
$$
Le us consider the second case for instance. For a fixed $\theta$ we have
\begin{eqnarray*}
\int_{0}^{r'_1}\module{h_n}^2e^{-2\reel\varphi}r^{2\beta-1}\module{\log r}^{k-2}\,dr&\leq &C(\theta) \int_{0}^{r'_1}\module{\chi f_n}^2e^{-2\reel\varphi}r^{2\beta+1}\module{\log r}^{k}\,dr
\end{eqnarray*}
with
\begin{eqnarray*}
C(\theta)&=&4\sup_{\rho\in[0,r'_1[}\int_{0}^{\rho}e^{-2\reel\varphi}r^{2\beta}\module{\log r}^{k-2}\,\frac{dr}{r}\cdot \int_{\rho}^{r'_1}\lefpar e^{-2\reel\varphi}r^{2\beta}\module{\log r}^{k}\rigpar^{-1} \frac{dr}{r}.
\end{eqnarray*}
As $\cos(\ell\theta-\tau)$ does not vanish in $[\theta_0,\theta_1]$, we may apply an argument similar to that of lemma \ref{lem:utile} to find $r_N(\varphi,\theta_0,\theta_1)$. We conclude that $C(\theta)$ is bounded by a constant independent of $\theta$, hence $\norme{h_n\otimes \wt e}\leq C\norme{f_ndr\otimes \wt e}$ and therefore $h_n\otimes \wt e$ tends to $h\otimes \wt e$ in $L^2(\Omega';L;d\vol)$.

At this step, we can repeat the proof for $\cH^{2}$ and thus get the proof of part (1) of lemma \ref{lem:rg1}. Let us end the proof of part (2).

We have
\begin{eqnarray*}
\frac{\partial h_n}{\partial \theta}&=&\chi g_n(r,\theta)-\int_{r}^{r'_1}\chi'(\rho)g_n(\rho,\theta)\,d\rho + \int_{r}^{r'_1}\chi\cdot\lefpar \frac{\partial f_n}{\partial \theta}-\frac{\partial g_n}{\partial \rho}\rigpar d\rho
\end{eqnarray*}
and one gets in the same way the convergence of ${\partial h_n}/{\partial \theta}$ in $L^2(\Omega;L;d\vol)$. This shows that $h\otimes \wt e$ and $dh\otimes \wt e$ are in $L^2(\Omega;L;d\vol)$.

We thus have $\omega-dh=\ell(\theta)d\theta$ with $\ell(\theta)d\theta\otimes \wt e\in L^2(\Omega;L;d\vol)$ and we are reduced to the case where $\omega=\ell(\theta)d\theta$. Then $\omega=0$ if $\cos(\ell\theta-\tau)>0$ on $]\theta''_0,\theta''_1[$, since $e^{-\reel\varphi}\rightarrow \infty$ exponentially when $r\rightarrow 0$, $\theta$ being fixed. We may thus assume that $\cos(\ell\theta-\tau)<0$ on $]\theta''_0,\theta''_1[$. Let $\theta_{1/2}\in{}]\theta_0,\theta_1[$ be such that $\sin(\tau-\ell\theta_{1/2})=0$, so $e^{-\reel\varphi}$ is increasing on $]\theta_0,\theta_{1/2}[$ and decreasing on $]\theta_{1/2},\theta_1[$.

Put then
\begin{eqnarray*}
u(\theta)&=&
\dpl\int_{\theta_{1/2}}^{\theta}\ell(t)dt.
\end{eqnarray*}
We have, using once more Hardy-type inequalities,
\begin{eqnarray*}
\int_{\theta_0}^{\theta_{1/2}}\module{u(\theta)}^2e^{-2\reel\varphi(re^{i\theta})}d\theta&\leq &C\int_{\theta_0}^{\theta_{1/2}}\module{\ell(\theta)}^2e^{-2\reel\varphi(re^{i\theta})}d\theta
\end{eqnarray*}
with
\begin{eqnarray*}
C&=&4\sup_{\theta\in[\theta_0,\theta_{1/2}]}\int_{\theta_{0}}^{\theta}e^{-2\reel\varphi(re^{it})}\,dt\cdot \int_{\theta}^{\theta_{1/2}}e^{2\reel\varphi(re^{it})}\,dt\\
&\leq &(\theta_{1/2}-\theta_{0})^{2}
\end{eqnarray*}
and analogously
\begin{eqnarray*}
\int_{\theta_{1/2}}^{\theta_{1}}\module{u(\theta)}^2e^{-2\reel\varphi(re^{i\theta})}d\theta&\leq &(\theta_{1/2}-\theta_{1})^{2}\int_{\theta_{1/2}}^{\theta_{1}}\module{\ell(\theta)}^2e^{-2\reel\varphi(re^{i\theta})}d\theta.
\end{eqnarray*}
We hence get
\begin{eqnarray*}
\int_{\theta_0}^{\theta_1}\module{u(\theta)}^2e^{-2\reel\varphi(re^{i\theta})}d\theta&\leq &C\int_{\theta_0}^{\theta_1}\module{\ell(\theta)}^2e^{-2\reel\varphi(re^{i\theta})}d\theta
\end{eqnarray*}
with $C$ independent of $r$. Therefore we have
\begin{eqnarray*}
\norme{u\otimes \wt e}^2&=&\iint\module{u}^2e^{-2\reel\varphi}r^{2\beta-1}\module{\log r}^{k-2}\,d\theta\,dr\\
&\leq& C \iint\module{\ell}^2e^{-2\reel\varphi}r^{2\beta-1}\module{\log r}^{k-2}\,d\theta\,dr\\
&\leq& C' \iint\module{\ell}^2e^{-2\reel\varphi}r^{2\beta-1}\module{\log r}^{k}\,d\theta\,dr\\
&=&C'\norme{\ell\otimes \wt e}^2.
\end{eqnarray*}
\end{proof}

\section{Proof of theorem \ref{th:poincare} (continued)}\label{sec:continued}
\subsection{}\label{subsec:good}
By the notation $\cM$, or $\wh\cM$, \etc., we now understand the $\cD$-module structure, so that the connection is included. We will use the decomposition $\wh\cM=\wh\cM^{\reg}\oplus \wh\cM^{\ir}$ (\cf.\ \T\ref{subsec:rappels}).

We have $\cH_{[D]}^{0}(\cM)=\cH_{[D]}^{0}(\wh{\cM})$ and we have seen that $\cM_{\min}$ is the kernel of the surjective morphism
$$
\cM\longrightarrow \wh{\cM}\longrightarrow \wh{\cT}=\cT.
$$
Moreover, we have $\wh{\cT}=\wh{\cT^{\reg}}$, where the last term is computed with $\wh\cM^{\reg}$. Thus $\DR(\cM_{\min})$ is isomorphic to the cone (shifted by $-1$) of the composed morphism (viewing $\DR(\wh{\cM})$ \etc. as a complex supported at the origin):
$$
\DR(\cM)\longrightarrow\DR(\wh{\cM})\longrightarrow \DR(\wh{\cT})\longrightarrow \DR(\wh{\cT^{\reg}}).
$$

Assume now (\S\T\ref{subsec:cassimple}-\ref{subsec:fingood}) that $\cM$ has a formal decomposition $\wh{\cM}\simeq\wh\cM^{\el}$ with $\cM^{\el}=\oplus_\varphi (\cE^{\varphi}\otimes \cR_\varphi)$. It will be easier to work with the metric $k_1$ instead of $k$ (the $L^2$ complexes are the same).

\subsection{}\label{subsec:cassimple}
For simplicity, we will first consider the case where $\cM$ has no regular component. We thus have $\cM_{\min}=\cM$ and $\cT=0$. We have to show that the natural morphism $\cL_{(2)}^\bbullet(\cM)\rightarrow \Db_{X}^{\bbullet}\otimes\cM$ is a quasi-isomorphism.

As explained in \T\ref{subsec:tilde}, it is enough to prove that the natural inclusion morphism
$$
\cL_{(2)}^{\bbullet}(\wt\cM)\HMRE{} \Db_{\wt\Delta}^{\mod0,\bbullet}\otimes_{\cA_{\wt\Delta}}\wt\cM
$$
is an isomorphism in $D^b(\CC_{\wt\Delta})$. This is a local problem on $\wt\Delta$, so we may assume, by the Hukuhara-Turrittin theorem recalled in \T\ref{subsec:lift}, that $\cM$ is elementary, as $\wt\cM$ and $\wt\cM^{\el}$ are locally isomorphic on $\wt\Delta$.

It is easy to see that this morphism induces an isomorphism on $\cH^0$. Moreover, it is known that $\cH^i(\Db_{\wt\Delta}^{\mod0,\bbullet}\otimes_{\cA_{\wt\Delta}}\wt\cM)=0$ for $i\geq 1$: indeed, this cohomology sheaf is equal to $\cH^i(\cA_{\wt\Delta}^{\mod0}\otimes_{\cA_{\wt\Delta}}\DR(\wt\cM))$, with $\cA_{\wt\Delta}^{\mod0}=\ker\lefcro \ov z\partial_{\ov z}:\Db_{\wt\Delta}^{\mod0}\rightarrow \Db_{\wt\Delta}^{\mod0}\rigcro$, and the latter group is zero, as indicated in \cite[p\ptbl211]{Malgrange91}.

Therefore, it is enough to show that $\cH^i\lefpar \cL_{(2)}^{\bbullet}(\wt\cM)\rigpar=0$ for $i\geq 1$. This follows from lemma \ref{lem:rg1} in \T\ref{sec:firstpart} by an easy extension argument, taking $\beta=0$.

\subsection{}\label{subsec:fingood}
Recall that (remark \ref{rem:zucker}-(2)), for regular connections, the theorem is proved by Zucker. Hence it follows from \T\ref{subsec:cassimple} that the theorem is true for elementary connections $\cM^{\el}$. It remains to show it for connections $\cM$ having a regular part in their elementary local model $\cM^{\el}$.

\begin{lemme}
Let $\theta^o\in S^1$. Any morphism $\wt p:\wt\cM_{\theta^{o}}\to\wt\cM_{\theta^{o}}^{\reg}$ lifting the projection $\wh\cM\to\wh\cM^{\reg}$ induces the same morphism
\[
\wt c:\cH^{1}(\cL_{(2)}^{\bbullet}(\wt\cM))_{\theta^{o}}\longrightarrow \cH^{1}(\cL_{(2)}^{\bbullet}(\wt\cM^{\reg}))_{\theta^{o}}.
\]
\end{lemme}

\begin{proof}
Let us fix a decomposition $\wt\cM_{\theta^{o}}=\wt\cM_{\theta^{o}}^{\reg}\oplus\wt\cM_{\theta^{o}}^{\ir}$ lifting the formal decomposition, and let $\wt p_{1}$ be the first projection. One may write $\wt p=\wt p^{\reg}+\wt p^{\ir}$ where $\wt p^{\reg}:\wt\cM_{\theta^{o}}^{\reg}\to \wt\cM_{\theta^{o}}^{\reg}$ is asymptotic to $\id$, and $\wt p^{\ir}:\wt\cM_{\theta^{o}}^{\ir}\to \wt\cM_{\theta^{o}}^{\reg}$ is asymptotic to $0$.

Then $\wt p^{\reg}$ is equal to $\id$, because $\wt p^{\reg}-\id$ is a horizontal section of the meromorphic connection $\cHom_\cO(\cM^{\reg},\cM^{\reg})_{\theta^o}^{{}^{\scriptstyle\sim}}$; as $\cHom_\cO(\cM^{\reg},\cM^{\reg})$ has a regular singularity, any horizontal section in a sector which is asymptotically equal to $0$ vanishes identically.

On the other hand, lemma \ref{lem:rg1} of \T\ref{sec:firstpart} shows that $\cH^1\lefpar \cL_{(2)}^{\bbullet}(\wt\cM^{\ir})\rigpar_{\theta^o}=0$.

Consequently, $\wt p$ and $\wt p_{1}$ induce the same morphism $\wt c$.
\end{proof}

This lemma implies that there exists a globally well-defined isomorphism
\[
\wt c: \cH^1\lefpar \cL_{(2)}^{\bbullet}(\wt\cM)\rigpar\isom \cH^1\lefpar \cL_{(2)}^{\bbullet}(\wt\cM^{\reg})\rigpar.
\]

Consider the following exact sequences of complexes:
$$
\xymatrix{0\ar[r]&\cL_{(2)}^{\bbullet}(\wt\cM)\ar[r]&\Db_{\wt\Delta}^{\mod0,\bbullet}\otimes \wt\cM \ar[r]^-{\wt a}&\wt\cC\ar[r]&0\\
0\ar[r]&\cL_{(2)}^{\bbullet}(\wt\cM^{\reg})\ar[r]&\Db_{\wt\Delta}^{\mod0,\bbullet}\otimes \wt\cM^{\reg} \ar[r]&\wt\cC^{\reg}\ar[r]&0
}$$
By Zucker we have $\bR e_*\wt\cC^{\reg}=\cT^{\reg}$. Moreover, $\Db_{\wt\Delta}^{\mod0,\bbullet}\otimes \wt\cM$ has cohomology in degree $0$ only and the morphism $\cH^0(\cL_{(2)}^{\bbullet}(\wt\cM))\rightarrow \cH^0(\Db_{\wt\Delta}^{\mod0,\bbullet}\otimes \wt\cM)$ is an isomorphism: indeed, this is a local statement on $\wt\Delta$ which is true for elementary local models as $\cM^{\el}$, so is also true for $\cM$. It follows that $\wt\cC$ has cohomology in degree $0$ only, this cohomology being isomorphic to $\cH^1(\cL_{(2)}^{\bbullet}(\wt\cM))$.

Consequently, the morphism $\wt c$ above induces an isomorphism $\wt c:\wt\cC\rightarrow \wt\cC^{\reg}$ in $D^b(\CC_{\wt\Delta})$.

We will now use the sheaves $\cA^{\leq 0}$ and $\wh\cO^{\rm Nils}$ defined on $S^{1}=e^{-1}(0)$ in \cite[p\ptbl61]{Malgrange91}, and denote the de~Rham complexes with coefficients in these sheaves by $\DR^{\leq 0}(\wt\cM)$ and $\DR^{\leq 0}(\wh\cM)$.

The natural morphism $\DR^{\leq 0}(\wt\cM)\rightarrow (\Db_{\wt\Delta}^{\mod0,\bbullet}\otimes \wt\cM)_{|S^{1}}$ is a quasi-isomorphism (\cf. \loccit.) because it induces an isomorphism of $\cH^0$ and the $\cH^i$ vanishes for $i\geq 1$.

The left part of the following diagram of complexes defines the morphism $\wt b$:
$$
\xymatrix{
\DR^{\leq 0}(\wt\cM)\ar[d]\ar[r]^-\sim\ar@{}[rrddd] |{\circlearrowright} &(\Db_{\wt\Delta}^{\mod0,\bbullet}\otimes \wt\cM)_{|S^{1}}\ar[r]^-{\wt a} \ar[rddd]^-{\wt b}&\wt\cC\ar[ddd]^{\wt c}_\wr\\
\DR^{\leq 0}(\wh\cM)\ar[d]&&\\
\DR^{\leq 0}(\wh\cM^{\reg})&&\\
\DR^{\leq 0}(\wt\cM^{\reg}) \ar[u]^-\wr\ar[r]^-\sim&(\Db_{\wt\Delta}^{\mod0,\bbullet}\otimes\wt\cM^{\reg})_{|S^{1}}\ar[r]&\wt\cC^{\reg}
}
$$

Up to a shift by $-1$, the complex $\DR(\cM_{\min})$ is identified with the cone of $\bR e_*\wt b$ and the complex $\cL_{(2)}^{\bbullet}(\cM)$ with the cone of $\bR e_*\wt a$.

To end the proof of the theorem under the assumption on the existence of a formal decomposition, it is therefore enough to show that $\wt b=\wt c\circ \wt a$. This is now a local problem on $S^1$ and we may use a local morphism $\wt\cM\rightarrow \wt\cM^{\reg}$ lifting the projection $\wh\cM\rightarrow \wh\cM^{\reg}$ to define $\wt c$ at the level of complexes and not only at the level of the derived category. With this definition of $\wt c$ the diagram clearly commutes.

\subsection{}
Let us now consider the general case. Let $\pi:\Delta_q\rightarrow \Delta$ be a cyclic covering of degree $q$ such that the meromorphic connection $\pi^+\cM$ has a formal decomposition. Fix the metric $k$ as in \T\ref{subsec:decform}. The theorem is thus true for $\pi^+\cM$.

First, one shows that $\cL_{(2)}^{i}(\cM)$ is a direct factor of $\pi_*\cL_{(2)}^{i}(\pi^+\cM)$. Indeed, one can write each section of $\pi_*\cL_{(2)}^{i}(\pi^+\cM)$ as a sum of terms $\omega\otimes \pi^*\wt e$, where $\wt e$ is a horizontal multivalued section of $\cM$ on $\Delta^*$ and where $\omega$ is expressed with the forms $dz/z$ and $d\ov z/\ov z$; the composition of the coefficients of $\omega$ with $\pi$ is seen to induce an inclusion $\cL_{(2)}^{i}(\cM)\subset \pi_*\cL_{(2)}^{i}(\pi^+\cM)$, if the metric on $\pi^+\cM_{|\Delta_q^*}$ is the inverse image of the metric $k_1$ on $\cM_{|\Delta^*}$ (it is comparable to the metric used in \T\ref{subsec:fingood}, hence induces the same $L^2$ spaces); the projection on $\cL_{(2)}^{i}(\cM)$ is given by the trace map $1/q\cdot\tr$.

The complex $\cL_{(2)}^{\bbullet}(\cM)$ is therefore a direct factor of $\pi_*\cL_{(2)}^{\bbullet}(\pi^+\cM)$. According to the theorem applied to $\pi^+\cM$, the complex $\cL_{(2)}^{\bbullet}(\pi^+\cM)$ is perverse on $\Delta_q$, being isomorphic to the perverse complex $\DR((\pi^+\cM)_{\min})$. The map $\pi$ being finite, the complex $\pi_*\cL_{(2)}^{\bbullet}(\pi^+\cM)$ is perverse on $\Delta$ as well as the direct factor $\cL_{(2)}^{\bbullet}(\cM)$.

Let $\zeta=e^{2i\pi/q}$ and let $\sigma$ be the automorphism of $\Delta_q$ sending $t$ to $\zeta t$. We have $\pi\circ \sigma=\pi$. It induces an automorphism $\wt\sigma$ of $\pi_*\cL_{(2)}^{\bbullet}(\pi^+\cM)$:
\begin{eqnarray*}
\pi_*\cL_{(2)}^{\bbullet}(\pi^+\cM)&=&(\pi\circ\sigma)_*\cL_{(2)}^{\bbullet}((\pi\circ\sigma)^+\cM)\\
&\simeq&\pi_*\sigma_*\cL_{(2)}^{\bbullet}(\sigma^+\pi^+\cM)\\
&\simeq&\pi_*(\sigma_*\sigma^{-1})\cL_{(2)}^{\bbullet}(\pi^+\cM)\\
&=&\pi_*\cL_{(2)}^{\bbullet}(\pi^+\cM).
\end{eqnarray*}
The perverse subsheaf $\cL_{(2)}^{\bbullet}(\cM)$ is then identified with $\ker(\wt\sigma-\id)$ (the kernel is taken in the perverse sense).

We may argue in a similar way with the terms appearing in the following diagram:
\[
\xymatrix{
\pi_*\cL_{(2)}^{\bbullet}(\pi^+\cM)\ar[r]&\pi_*(\Db_{\Delta_q}^{\mod0}\otimes \pi^+\cM)\ar[r]&\pi_*\cC(\pi^+\cM)\ar[d]^-{c}\\
\pi_*\DR((\pi^+\cM)_{\min})\ar[r]&\pi_*\DR(\pi^+\cM)\ar[r]\ar[u]^-{\wr}&\pi_*\DR(\cT_{\pi^+\cM})
}
\]
The theorem for $\cM$ will follow from the result for $\pi^+\cM$ and from the fact that $\pi_*c_{\pi^+\cM}$ is compatible with $\wt\sigma$. Let us show this last statement: the decomposition $\pi^+\wh\cM=(\pi^+\wh\cM)^{\reg}\oplus (\pi^+\wh\cM)^{\ir}$ is invariant under the action of $\sigma$; taking the mean of any local liftings $\wt{\pi^{+}\cM}\to\wt{\pi^{+}\cM^{\reg}}$ of the projection $\pi^+\wh\cM\rightarrow (\pi^+\wh\cM)^{\reg}$ at the points $\theta^o+2ik\pi/q$ gives thus another such local lifting, which is invariant under the action of $\sigma$; as $\wt c$ does not depend on the choice of the local lifting, we conclude that $\pi_*c_{\pi^+\cM}$ is compatible with $\wt\sigma$.\hfill\qed

\subsection*{Acknowledgements}
Part of this work has been achieved during a visit of the author at the Tata Institute (Bombay) in the framework of the Indo-French exchange program IFCPAR. The author thanks Nitin Nitsure and Olivier Biquard for useful discussions.

\providecommand{\bysame}{\leavevmode\hbox to3em{\hrulefill}\thinspace}

\end{document}